\documentstyle[11pt]{article}
\onecolumn
\newtheorem{lem}{{\sc Lemma} }[section]
\newtheorem{thm}{{\sc Theorem} }[section]
\newtheorem{cor}{{\sc Corollary} }[section]

\title{\LARGE\bf
   Finite Semigroups of Constant Rank, \\[1ex]
   and the five Basic State Machine types}
\author{ {\sc Nico F. Benschop} \\[1ex]
  {\it n.benschop@chello.nl - Amspade Research} - Geldrop, ~The Netherlands}
\date{ {\it IFIP - TC10 Workshop on "Logic and Architecture Synthesis"\\
       30 May - 1 June 1990, Paris (publ: North-Holland, p167-176} ) }

\begin{document}
\maketitle

\begin{abstract}
Constant Rank ($CR$) state machines play an important role in the general
structure theory of Finite State Machines. A machine is of constant rank if
each input and input-sequence maps the state set onto the same number of
next states. $CR$-machines are analysed via their sequential closure
(semigroup), which is a {\bf simple semigroup :} a semi-
direct product $(L \times R) ~* G$ of a left- and a right-copy semigroup,
and a group.\\
So in general a $CR$-machine is a composition of: ~a {\bf branch-},
a {\bf reset-} and a {\bf permutation} machine, ~which are three of the
five basic types of state machines [1].
~(Original title:~ "The Structure of Constant Rank State Machines")
\end{abstract}

\section{Introduction: ~~Sequential closure and rank}  % sect 1

A brief review of [1] is necessary to set up the required concepts. A state
machine $M(Q,A)$ with stateset $Q$ and input alphabet $A$ is a function
~$M: Q \times A ~\rightarrow~ Q$, which maps present state and input to next state.
It is specified by a state transition table with $|A|$ columns and $|Q|$ rows.
Each input $a \in A$~ is interpreted as a function ~$a: Q \rightarrow Q$, mapping
stateset $Q$ into itself, called a {\bf state transform}, or in short: a transform.

Sequential composition $ab$ of two transforms $a$ and $b$ is defined by ~$q(ab)=(qa)b$,
for all $q \in Q$. In other words, in state $q$ first apply input $a$ to get state
$qa$, then apply $b$~ which yields state $(qa)b = q(ab) = qab$. Notice the
left-to-right notation of this {\bf function composition}, with stateset $Q$ as
domain and codomain. Two input sequences over $A$ are defined {\bf equivalent} if
they yield the same $Q$-transform: ~$a=b$~ iff ~$qa=qb$~ for all $q \in Q$.

The {\bf sequential closure of $M$}, called {\bf semigroup $S$}, is the (finite)
set of $Q$-transforms generated by all sequences over $A$, denoted ~$S=A^+/Q$.
Here ~$A^+$~ denotes the infinite semigroup of non- empty strings, length
$\geq$ 1 over alphabet $A$, under string concatenation.

Closure $S$ of machine $M$ is a finite semigroup (of order $|S| \leq n^{n}$,
if $M$ has $n$ states) since transform composition is {\bf associative}:
$a(bc)=(ab)c$ for all $a,b,c \in S$, which is clear from above definition of
transform composition. Input-strings with the same $Q$-transform are defined
{\it equivalent} with respect to machine $M$, so the transform representation
of each element of $S$ is unique. State transform $x: Q \rightarrow Q$~ is a
function defined on state set $Q$, which is both domain and co-domain.
To state transform $x$ correspond:

- {\bf range ~$Qx$} ~which is the set of function values (next states), and

- {\bf partition $Px$} ~equivalences states that map onto the same next state.

- {\bf rank $r(x)$} = the order $|Qx|$ of its range = the number of partition
blocks.

\begin{lem} : ~~(non-increasing rank property) % lem 1.1

(a) ~Left composition ($x.$) does not increase range :
      ~~$Qxy \subseteq Qy$~~( $\subset$ : subset of)

(b) ~Right composition ($.y$) does not refine partition:
       ~$Pxy \geq Px$~~( $>$ : coarser than)

(c) ~Rank does not increase under transform composition:
        ~$r(xy) \leq r(x)$~ and ~$r(xy) \leq r(y)$

(d) ~All elements $x$ with $rank(x) \leq k$ form a subsemigroup
        which is an ideal ~$Z_k$ ~of~ $S$.
\end{lem}
\begin{proof}
(a) ~$Qxy \subseteq Qy$~ follows from set inclusion and associativity.

     ~$Qx \subseteq Q$ ~for all $x$, and right composition with $y$ yields:
   $(Qx)y=Q(xy) \subseteq Qy$.

(b) ~$Pxy \geq Px$~ follows from associativity and right composition of states
$i,j$ that are equivalent under~ $.x:$~ ~$ix=jx$ implies ~$ixy=jxy$~ for all $y$.
~So ~$i \equiv_x j$~ implies ~$i \equiv_{xy} j$.

(c) ~This {\bf monotone rank property} follows directly from (a) and (b),

because range ordering (a) implies rank ordering $|Q(xy)| \leq |Qy|$,
	so ~$r(xy) \leq r(y)$,

~and partition ordering (b) implies rank ordering $|P(xy)| \geq |Px|$,
	so ~$r(xy) \leq r(x)$.

(d) ~It follows immediately that if $x$ and $y$ have rank $\leq$ k, then so does
composition $xy$. This closure property means that all elements of rank not
exceeding k form a subsemigroup $Z$ of $S$. In fact, composition of any element
$z \in Z$~ with any element $s \in S$~ yields ~$zs$~ with ~$r(zs) \leq r(z) \leq k$,
sothat $zs \in Z$. The same holds for $sz$. Hence $Z$ is both left- and right ideal,
that is an ideal of $S$ ~with~ $ZS \subseteq Z$~ and ~$SZ \subseteq Z$~
(see def-2 next section). \end{proof}

Basically, this paper tries to render results from semigroup structure and their
state representation better accessible for state machine decomposition purposes.
In fact, the earliest known result in semigroup theory (Suschkewitch, 1928 [2,
p207]) is on the structure of the minimal ideal of a semigroup, essentially
our theorem 4.1.

\section{Basic machines and simple semigroups} % sect.2

Machine decomposition is seen as implementing a machine as a network of
$smaller$ machines. Semigroups, as the sequential closures of state machines,
are essential for the {\bf equivalencing} and {\bf ordering} of machines.
Two machines are defined to be
equivalent if they have isomorphic semigroups. Two machines are ordered
~$M_1 \leq  M$ if their closures are ordered ~$S_1 \leq S$, meaning that $S_1$
is (isomorphic to) a subsemigroup of $S$.\\
{\bf Def 1}: a minimal or {\bf basic machine} has a closure with
{\it no proper subsemigroup}.

In [1] it is shown
that the minimal number of generators, the 'dimension' $dim(S)$, of a
$basic$ semigroup $S$~ is either one ($iterative$ structure $S=a^*/Q$),
or two ($idempotent$ generators $S=\{a,b\}^*/Q$ with $a^2=a, ~b^2=b$).
Because if at least three generators were required, any two of them would
generate a proper subsemigroup. And  if two are required, then no generator
can generate more than itself ($idempotent$ or 'invariant') since otherwise
$|a^*|>1$~ yields a proper (iterative) subsemigroup. Such idempotent
pair can generate either a commutative basic $S=H_2$ of two ordered
invariants, or one of two non-commutative basic left- or right- copy
semigroups $L_2$ or $R_2$, with $ab=a$ resp. $ab=b$. Iterative $S=a^*/Q$
are $basic$ if they are $periodic$ (see section 3) and of prime order
($C_p$), or $monotone$ ~(type $U$) of order 2.

The {\bf five basic state machines} with semigroups of order two, are derived
in [1], with their interpretation as the elementary digital functions of type:
{\it logic ($H$), arithmetic ($C,U,H$)} and $memory$ ($L,R$). A semigroup $S$
is also a state machine ~$M(S,S)$~ with itself as inputset and state set. For
unique representation by state transforms (distinct columns), one extra state
suffices if some columns are equal in the $S \times S$ composition table,
see tables $U_2$ and $L_2$. Components $C_2$~ and ~$U_2$~ have a single
generator '1', the others have two invariant generators ~$a^2=a$.

{\bf Def 2:} a semigroup is of {\bf constant rank} ($CR$) if it can be represented
by transforms of equal rank. A state machine is of constant rank if its closure
is a $CR$-semigroup.

Three basic components are of constant rank, namely $L_2,~R_2$~ and $C_2$.
They are the smallest cases of the following three types of {\bf constant
rank semigroups}:

$L:$ ~~Left-copy semigroup ~~with ~$ab=a, ba=b$ for all $a,b \in S$
   ~($n$-branch, $n$+1 states)

$R:$  ~~Right-copy semigroup ~~with ~$ab=b, ba=a$ for all $a,b \in S ~(n$-reset,
    $n$ states)

$G:$ ~~Group (permutation machine: permutes ~$n$~ states,  ~$|G|\leq n!$ )

All three are special cases of the following general type of semigroup [2, p5]:

{\bf Def 3:} ~an {\bf ideal} of a semigroup $S$ is a subset $Z$ with
~$SZ \subseteq Z$~ and ~$ZS \subseteq Z$. \\ \hspace*{1.5cm}
A semigroup is called {\bf simple} if it has {\it no proper ideal.}

An ideal is like a multiplicative 'zero' ($a$.0=0 for all $a$) or 'trap'.
Notice that $U_2$ (monotone counter with a final state) and $H_2$ (hierarchy
of two ordered invariants, see next section) are not simple semigroups, nor
are they of constant rank. In general they model the {\it monotone sequential}
aspects and {\it combinational logic} aspects of state machines respectively.

\begin{verbatim}
C2| 1 0       U2| 1 0      H2| 1 0     L2| 1 0     R2| 1 0
--+----       --+----      --+----     --+----     --+----    Closure
 1| 0 1        1| 0 0       1| 1 0      1| 1 1      1| 1 0     Tables
 0| 1 0        0| 0 0       0| 0 0      0| 0 0      0| 1 0
               2| 1 0                   2| 1 0                 Fig.1a

  .<-.                       o 1        /->o 0
 o-->o       o-->o-->o       :       2 o            o    o    State-
 1   0       2   1   0       o 0        \->o 1      0    1    Diagrams

2-counter    2-counter    AND, OR     2-branch   set/reset    Component
 periodic    monotone     isomorph      mux        D-FF       Functions
add(mod 2)   converge    mpy(mod 2)   if-else    assign :=

<---- iterative a* ----> <---- invariant : aa=a -------->     Algebraic
                         <- LOGIC ->  <-select-> <-store->    Properties
<--- ARITHMETIC : commutative ----->  <- MEMORY non cmt ->
                                                               Fig.1b
\end{verbatim}

\begin{cor} % cor 2.1
~~A {\bf simple} semigroup is of {\bf constant rank}. \end{cor}

This follows directly from lemma 1.1d, since otherwise the elements of
minimum rank would form a proper ideal. In fact, it will be shown that any
simple semigroup is a semi-direct product $(L \times R) ~*~ G$ of the three
basic types of simple semigroups $L,~R,~G$.

So a general $CR$-machine is the parallel composition of a {\it branch}
machine, a {\it reset} machine and a {\it permutation} machine. In a way,
this is a conservation law of sequential logic.

\section{Iterations: ~~monotone, ~periodic, ~invariant} % sect.3

Iteration in a semigroup $S$ is the repetition $a^i$ of a single element. By
virtue of associativity, the result is a unique element in $S$, independent
of bracketing. The closure of a single element $a \in S$ is the finite set
of its {\bf iterations} $a^+ = \{a^i, ~i=1..n\}$ which in general has a
tail-cycle structure ( $/Q$ is omitted if no confusion can arise):

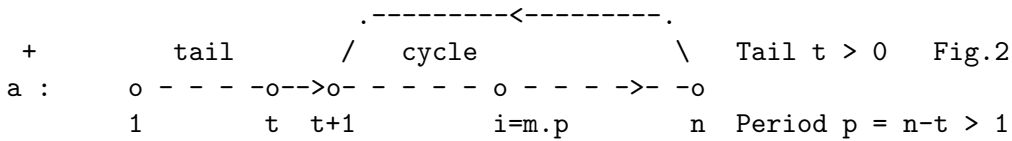
\begin{verbatim}
                       .---------<---------.
 +         tail       /   cycle             \   Tail t > 0   Fig.2
a :     o - - - -o-->o- - - - - o - - - ->- -o
        1        t  t+1         i=m.p        n  Period p = n-t > 1
\end{verbatim}

Since $a^+$ is finite, there is a smallest n for which $a^{n+1} = a^{t+1}$
with {\bf tail}$(a)=t, ~0 \leq t<n$~ and {\bf period}$(a)= p = n-t$. There is
precisely one {\bf invariant} $a^{i} = (a^{i})^{2}$ where $i=mp$ is the first
and only multiple of $p$ in the cycle, and $a^{k}=a^{k+p}$ for $k>t$.

An element of semigroup $S$ is called {\bf periodic [monotone]} if its closure
has no tail, $t=0$  ~[ no cycle, $p=1$ ]. Clearly, {\it invariants} $aa=a$
are the only elements which have both properties. Elements which have a tail
and a cycle are called $aperiodic$.

{\bf Def 4}: ~a pair $e,z$ of {\bf commuting invariants}: $ez=ze, e^2=e,
z^2=z$, is said to be {\bf ordered} $e \geq z$ ~when $ez=ze=z$~ hence~ $e$ is
left- and right- identity for $z$. This relation is easily seen to be reflexive,
anti-symmetric and transitive [2, p23], so a partial ordering.

\subsection{Ordered Invariants: ~~H} % sect 3.1

It will be shown that any {\it simple semigroup} $S$, being of constant rank,
contains only $periodic$ elements. Moreover, its invariants are not ordered
but are all $equivalent$ in some sense. So basic components of type $U_2$
(monotone iteration)and $H_2$ (hierarchy of ordered invariants, or combinational
logic) do not occur. In fact it turns out that $S$ is a disjoint union of
isomorphic groups $G$, with identities forming a direct product of a left-copy
$L$ and a right-copy $R$ semigroup.

\begin{lem} : ~~(ordered invariants)\\  % lem 3.1
The ordering of commuting invariants $z \leq e$ is their range ordering:~
     $Qz \subseteq Qe$, \\ \hspace*{2cm} hence:
-- distinct commuting invariants have distinct ranges, and \\ \hspace*{3cm}
-- ordered invariants $z<e$ have ordered ranks $r(z) < r(e)$.
\end{lem}
\begin{proof}
Let invariants $z$ and $e$ be ordered $z \leq e$, then $e$ is identity for
$z: ez=ze=z$, so their ranges are ordered because $Qz=Q(ze)=(Qz)e \subseteq Qe$.
Notice that $ze=z$ suffices: $e$ is right identity for $z$. Conversely, for commuting
invariants: $Qz \subseteq Qe$ implies $z \leq e$. This follows from the {\it state
transform structure of an invariant}~ $e: qee=qe$ means that each state $q$ maps to
a state $qe$ which is fixed under $e$. In other words, no state chains of length
$>1$ occur in the state transition diagram of $e$.

Range $Qe$ is the set of {\bf fixed states} of $e$. Now, if $Qz \subseteq Qe$
then $z$ maps each state $q$ into a fixed state of $e: (qz)e=qz$ for all $q$,
so $ze=z$. Since by assumption $e$ and $z$ commute, we have $ez=ze=z$, which
means $z \leq e$. Clearly, if $Qe=Qz$ for commuting invariants $e$ and $z$,
then $e \leq z$~ and ~$z \leq e$, and hence $e=z:$ commuting invariants with
the same range are equal. \end{proof}

\begin{cor} ( anti-commutative )  % cor 3.1

~~~~~A {\bf simple} semigroup $S$ has {\it no ordered invariants},
      ~and no pair of invariants commutes.
\end{cor}
\begin{proof}
Ordered invariants have different ranks according to the previous lemma.
Let $k$ be the lowest rank of an ordered pair of invariants. Then, with
lemma 1.1d, $S$ has a proper ideal consisting of all elements with
$rank \leq k$,~ which contradicts $S$ being simple.

If invariants $e,f$ commute: $ef=fe$, then their composition $d=ef$ is also
invariant: $d^{2}=d$ since $ef.ef= ef.fe= e.ff.e= e.fe= e.ef= ee.f= ef$.
Moreover: $d$ is ordered under $e$, since $ed= eef= efe= de=d$~ so ~$d \leq e$,
~and similarly $d \leq f$.

It is easily verified [2, p24] that~ $d$ is the {\bf greatest lower bound}
~or~ {\bf meet} of $e$ and $f$.  So a commuting pair of invariants is either
ordered, or their composition is ordered under both, contradicting simple $S$.
Hence no pair of invariants commutes.
\end{proof}

So a semigroup of {\it commuting invariants} is partially ordered set where
each pair has a meet (set intersection), called a lower semilattice,
with a global zero. For $n$ states, there are at most $2^n$ commuting
invariants (Boolean lattice).

\subsection{Equivalent Invariants: ~~L, R} % sect 3.2

Consider now the invariants of a simple semigroup S. They do not commute (cor.3.1).
Invariants that do not commute may be equivalent in the following sense:

{\bf Def 5: ~~Equivalent Invariants} \\
-- ~Invariants $a,b$ forming a left- [right-] copy semigroup $L_2 ~[R_2]$\\
    are {\bf left- [right] equivalent}, written $aLb$ ~[$aRb$]

-- ~Invariants $a,b$ are {\bf ~equivalent}, ~denoted $a \sim b$,
    ~~if they are left- or right equivalent: \\ 
  either $directly$, forming $L_2$ or $R_2$,~
  or $indirectly$: alternating $L$ and $R$ via other invariants.

\begin{lem} % lem 3.2
Consider invariants $a,b$ in any semigroup $S$, represented over stateset $Q:$

(a) ~Equivalent invariants have equal rank: $a \sim b \Rightarrow |Qa|=|Qb|$,
\\ \hspace*{.5cm} but equal rank is not sufficient for equivalence: see (b)

(b) Let $(ab)^k=ab$ and $(ba)^k=ba$, with invariants $(ab)^{k-1}=ab^0$
      and $(ba)^{k-1}=ba^0$, \\ \hspace*{1cm} with max-subgroups
      $G_{ab^0}=\{x^i=ab^0$ ~for some ~$i>0\}$ resp. $G_{ba^0}$, then:

~~~ if $k$=2: $\{a,b,ab,ba\}$ are 2 or 4 invariants of equal rank forming
         $L_2, ~R_2$~ or ~$L_2 \times R_2$,

~~~ if $k>$2 this structure holds for max-subgroups
       $\{G_a,G_b,G_{ab^0},G_{ab^0}\}$ under set product.
\end{lem}
\begin{proof}
{\bf (a)} ~There are three cases of equivalence for invariants $a,b:$
left-, right- and indirect equivalence. In the first two cases of "direct"
equivalence, rank-lemma 1.1 yields: \\ \hspace*{.5cm}
 $aLb$ ~implies ~$r(a)=r(ab) \leq r(b)$~ and ~$r(b)=r(ba) \leq r(a)$,
~sothat ~$r(a)=r(b)$; \\ \hspace*{.5cm}
$aRb$~ implies ~$r(a)=r(ba) \leq r(b)$~ and ~$r(b)=r(ab) \leq r(a)$,
~sothat ~$r(a)=r(b)$.

Hence left- or right equivalent invariants have the same rank.
Transitivity holds in both cases. For instance let $aLx ~(ax=a, xa=x)$
and ~$xLb ~(bx=b, xb=x)$~ then $aLb$, ~since ~$ab= ax.b= a.xb= ax= a$~,
and similarly ~$ba=b$. Also right equivalence is transitive.

If $aLc$ and $cRb$, where $c$ differs from $a$ and $b$, ~then $a,b$
are not directly left- or right equivalent, yet they are $indirectly$
equivalent, denoted ~$aLRb$. ~Here ~$LR$~ is an equivalence relation,
easily verified to be reflexive, symmetric and transitive.
If $a$ and $b$ are indirectly equivalent, via other invariants, 
then they have the same rank by transitivity.

{\bf (b)} ~There are several cases: direct and indirect equivalence,
with either $k$=2 or $k>2$.

For $k$=2, in the $direct$ equivalent case ~$aLb$~ and ~$aRb$~ the elements
~$ab$~ and ~$ba$~ are not different from $a$ and $b$, ~forming ~$L_2$~ and
~$R_2$~ respectively.
For $indirect$ equivalence of invariants $a$ and $b$, and in case $k$=2 the only
other intermediate elements are invariants $ab$ and $ba$, with $aba=a$ and
$bab=b$, ~seen as follows. Invariants $a,b$ must have equal rank: $|Qa|=|Qb|$
(lemma 3.2a), hence exact equality holds in $(Qa)b \subseteq Qb$, ~so ~$Qa.b=Qb$
(*) and similarly ~$Qb.a=Qa$ (**). Composing both sides of (*) on the right by
$.a$ and applying (**) yields ~$Qa.ba= Qba= Qa$. ~So sequence $.ba$~ permutes
~$Qa \rightarrow Qa$. \\
Since ~$ba$~ is invariant, this is the identity permutation, hence $(qa)ba=qa$~
for all $q$, meaning ~$aba=a$. Similarly, invariance of ~$ab$~ implies ~$bab=b$.

So strings of $length>2$ are equivalent to strings of $length \leq 2$,
which are just ~$a,b,ab,ba$, ~forming a closure of four invariants, with the
next equivalences (using ~$aba=a, bab=b):$

-- $aRab$~ since $a.ab=aa.b=ab$ and ~$ab.a=a, ~~abLb$~ since $ab.b=a.bb=ab$ and $b.ab=b$,

-- $bRba$~ since $b.ba=bb.a=ba$ and ~$ba.b=b, ~~baLa$~ since $ba.a=b.aa=ba$ and $a.ba=a$.

These relations are depicted in a rectangular form in the figure 3.
The four elements ~$\{a,b,ab,ba\}$~ form an invariant semigroup with direct
product structure $L_2 \times R_2$.

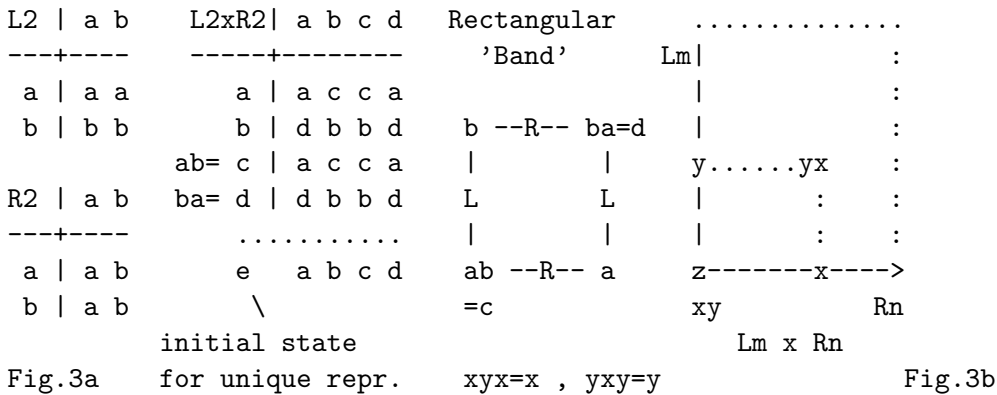
\begin{verbatim}
L2 | a b    L2xR2| a b c d   Rectangular     ..............
---+----    -----+--------     'Band'      Lm|            :
 a | a a       a | a c c a                   |            :
 b | b b       b | d b b d    b --R-- ba=d   |            :
           ab= c | a c c a    |        |     y......yx    :
R2 | a b   ba= d | d b b d    L        L     |       :    :
---+----       ...........    |        |     |       :    :
 a | a b       e   a b c d    ab --R-- a     z-------x---->
 b | a b        \             =c             xy          Rn
          initial state                         Lm x Rn
Fig.3a    for unique repr.    xyx=x , yxy=y                Fig.3b

Image = S / congruence : L2 = S/{a=c,b=d}; R2 = S/{a=d,b=c}
\end{verbatim}

$L_2 \times R_2$ is represented by a two-component code:~ $x=[x_1,x_2],
~y=[y_1,y_2]$ \\with~ $xy=[x_1,y_2]$ and $yx=[y_1,x_2]$.

In other words, the direct product $L_2 \times R_2$~ (for $k$=2) follows from
two complementary congruences (preserved partitions), illustrated by figure 3.
Denote ~$ab=c$~ and ~$ba=d$, then~ $\{a=c,b=d\}$~ with image $L_2$, and
~$\{a=d,b=c\}$~ with image $R_2$. The direct product is implemented by two
independent components~ $x=[x1,x2]:$~ the first composes as $L_2$ and the
second as $R_2$.

The left- and right equivalences can be plotted pairwise in the plane as
shown in fig 3, which also gives the composition tables of $L_2, ~R_2$~ and
~$L_2 \times R_2= \{a,b,ab,ba\}$. From this rectangular display follows the
term {\bf diagonal equivalence} for two indirectly equivalent invariants,
since this is the only other form of equivalence. It is denoted by ~$xDy$~
where $x$ and $y$ are obtained by commutation: $x=ab$~ and ~$y=ba$~ for some
$a$ and $b$, themselves being diagonal equivalent $aDb$, with
~$a=aba=abba=xy$~ and ~$b=bab=baab=yx$. Diagonal equivalence occurs in pairs:
if ~$aDb$~ then ~$ab D ba$, and vice versa.

The above analysis for $k$=2 can be generalized simply to $Lm \times Rn$
for $m.n$ invariants, with each invariant pair forming either $L_2$ or $R_2$
or $L_2 \times R_2$.

If $k>$2 in $(ab)^k=ab$ and $(ba)^k=ba$, then $ab$ and $ba$ are not invariant,
generating invariants $(ab)^{k-1}=ab^0$ and $(ba)^{k-1}=ba^0$ in a $k$-1 cycle,
with $(aba)^k=a$ and $(bab)^k=b$. The resulting structure is in general a
semi-direct product $(Lm \times Rn)* G$~ with a group $G$ as subgroup of $S$,
occurring $m.n$ times, to be derived next. In case $G$ is also an image of $S$,
then $S$ is direct product $(Lm \times Rn) \times G$.

Without going into much detail [2, Vol.I, appx]: each idempotent $a \in S$,
interpreted as left- or right- multiplier, yields (principle) subsemigroups
$aS$~ and ~$Sa$, respectively represented in the composition table of $S$ by
the rows and columns (fig.3). Each invariant $a$ is the identity of a maximal
subgroup $G_a = aSa$, the intersection of $aS$~ and ~$Sa$, while ~$aSb$~
contains ~$ab$~ and its invariant $(ab)^{k-1}$ as max-subgroup identity.
One readily verifies that all max-subgroups are isomorphic. Equivalencing
each to one congruence part, with $G_{ab}=G_a~G_b$, yields image ~$Lm \times Rn$~
where $m$ and $n$ represent the number of  max-subgroups in $S$ forming left-
resp. right- copy semigroups $Lm$ and $Rn$ as image. Notice that if the product
of invariants is not invariant, $Lm \times Rn$~ is not a sub-semigroup of $S$.

On the other hand: although $G$ occurs $m.n$~ as subgroup of ~$S$, it may
also be an image group~ $G=S / (Lm \times Rn)$, with $S=(Lm \times Rn) \times G$
as direct product. If this is not the case, so $G$ occurs as subgroups but not as
image of $S$, then $S$ is said to be a semi-direct product $S = (Lm \times Rn)~*~G$.
\end{proof}

The table of~ $L_2 \times R_2$~ viewed as state machine has two pairs of equal
columns (inputs $a=d$ and $b=c$), so an extra initial state $e$ is needed for
a unique state transform representation.

\begin{lem}:\\ % lem 3.3
{\bf (a)} ~~In any {\bf invariant} semigroup $S: ~~~a  \geq aba$.

{\bf (b)} ~~$a>aba$ ~for some ~$a,b$~ only if $S$ is not of constant rank, so~
\\.~~~~~~$a=aba$~ for all $a,b$ ~~iff~~ $S$ is invariant of constant rank.
\end{lem}
\begin{proof}
(a) ~We need to show that ~$a$~ commutes with $aba$, and is left- and right
identity for ~$aba$. Both follow directly from $aa=a$~ and ~$a.aba= aba =aba.a$.

(b) ~If $S$ is not of constant rank, then the minimum rank invariants form a
proper ideal $Z \subset S$ (lemma 1.1d), and there is an ordered and commuting pair
of invariants. Consider invariants ~$a \in S-Z$~ and ~$b \in Z$, then invariant
~$aba$~ is also in ~$Z$~ and has the same (minimal) rank as $b$, so
$rank(a) > rank(aba) = rank(b)$. Hence strict ordering ~$a>aba$~ holds.
\end{proof}

The rectangle of equivalent pairs of invariants generalizes to ~$Lm \times Rn$,
with ~$m,n \geq 2$. The $mn$ invariants form an  ~$m \times n$~ matrix, where
~$L$- ~[$R$-]~ equivalence holds between elements in the same column [row].
This is the general structure of a constant rank invariant semigroup (also
called a rectangular 'band'):

\begin{thm} % thm 3.1
~~The following conditions on a finite semigroup $S$ are equivalent:

(a) ~$S$ is {\it anti commutative} ~(no two elements commute: $ab=ba$ implies $a=b$).

(b) ~$S$ is {\it invariant} and of {\it constant rank}.

(c) ~$aba=a$ for all $a,b$ in $S$.

(d) ~Each pair $a,b$ of invariants in $S$ is equivalent: either directly,
    forming ~$L_2$~ or ~$R_2$, \\ \hspace*{.4in}
    or indirectly (diagonal) via $ab$ and $ba$ forming ~$L_2 \times R_2$.

(e) ~$S$~ is a direct product ~$Lm \times Rn$~ of a left- and a right copy
     semigroup~ $(m,n \geq 1)$.
\end{thm}
\begin{proof}
(a) $\Rightarrow$ (b) : an anti- commutative semigroup $S$ is invariant,
because any iteration class~ $x^+$ ~is a commutative subsemigroup, so~
$|x^+|=1$~ for all $x$, so each element of $S$ is invariant.
Moreover, $S$ is of constant rank; otherwise some pair of invariants $a,b$
would be properly ordered (lemma 3.3b) and thus commute, contradicting $S$
being anti-commutative.

(b) $\Rightarrow$ (c) : lemma 3.3b.

(c) $\Rightarrow$ (d) : $aba=a$ for all ~$a,b ~\rightarrow$ ~pairwise
   $L$-, $R$- ~or~ $D$- equivalent (lemma 3.2b).

(d) $\Rightarrow$ (e) : Pairwise equivalence in $S$ implies the direct
product structure $Lm \times Rn$~ with $m,n \geq 1$~ as follows. If $S$
contains only left- equivalent invariants then ~$S=Lm$~ where $m=|S|$
~and~ $n=1$. The other trivial case occurs when $S$ contains n right
equivalent invariants, and no left equivalence holds: $S=Rn$ with $m$=1
~and~ $n=|S|$.

If both left- and right equivalences occur, the ~$Lm \times Rn$
rectangular structure (fig.3b) is seen as follows. Take any invariant
$z$ and form two subsets: $Lz$ with all elements $y$ that are left
equivalent $yLz$ to $z$, and $Rz$ containing all $x$ with $xRz:$ right
equivalent to $z$. They intersect only in $z$, because if $w$ is left-
and right equivalent to $z$, then it cannot differ from $z: w=wz=z$.
~$Lz$~ and ~$Rz$~ are left- and right copy subsemigroups of $S$. Let
the orders be respectively $|Lz|=m$~ and ~$|Rz|=n$. Pairwise
equivalence implies n copies of $Lz$ which form a congruence $\lambda$
of $S$ with image $S/ \lambda= Rn$. Similarly, congruence $\rho$
consists of $m$ copies of $Rz$, yielding image $S/ \rho =Lm$. Since
no pair of invariants can be both left- and right equivalent,
congruences $\lambda$ ~and~ $\rho$ are orthogonal: ~$S= Lm \times Rn$.

(e) $\Rightarrow$ (a) : semigroup $S= Lm \times Rn$ consists of
pairwise equivalent invariants. Then it is anti- commutative which
means that no pair commutes. For assume that one pair of distinct
invariants $a,b$ commutes: $ab=ba$, then they are either ordered
$a<b$~ or ~$a>b$ (in case $ab$ is $a$~ or ~$b$), or their product
is a third invariant $c=ab=ba$, their {\bf meet}, that is ordered
$c<a$~ and ~$c<b$. Either case contradicts pairwise equivalence.
\end{proof}

Notice that rather general conditions (a)(b) imply a very regular
structure (e), which is due to the strong properties of $finite$
(rank) $associative$ (semigroup) algebra.

\section{Maximal Subgroups: ~~periodic $G$}  %sect 4

\begin{lem}  % lem 4.1
~~For the {\bf iterations} ~$a^i$~ of a semigroup element $a$ with increasing $i:$

-- the tail elements (if any) reduce strictly in rank, ~~~and

-- the cycle elements (at least one: the invariant of $a$) have
   constant minimum rank.
\end{lem}
\begin{proof}
Consider the successive ranges ~$Qa^i$~ which, due to range lemma 1.1a,
form a reducing inclusion chain of subsets of $Q$. Each range is
contained properly in the previous one until the cycle is reached
at $i=t+1$. As soon as two successive ranges are equal, then so are
all next ranges: $Qa^i=Qa^{i+1} ~\rightarrow~ Qa^{i+1}=Qa^{i+2}$, etc.
(compose left and right by $a$). Once the cycle is reached, the minimum
rank is obtained: the initial tail ranks decrease strictly, and all
periodic elements in thecycle have equal and minimal rank. \end{proof}

\begin{cor}  % cor 4.1
~In a simple semigroup $S$ every element is periodic (has no tail).
\end{cor}

This follows directly from the previous lemma and lemma 1.1d, because if an
element of S had a tail, then its iterations would have different ranks, which
contradicts the constant rank property of a simple semigroup. ~~~To show that a
simple semigroup is a disjoint union of isomorphic groups, we first need:

\begin{lem} % lem 4.2
~~( Maximal subgroups ) ~Let $S$ be a semigroup, then:

(a) ~Periodic elements generating the same invariant $e$ form a maximal
	subgroup of $S$, \\ \hspace*{1cm} ~called the group G$_{e}$ ~on~ $e$.

(b) ~Equivalent invariants ~$a \sim b$~ have isomorphic groups ~$G_a \cong G_b$:

~~~~~~ if ~$aLb$~ via isomorphism ~~$a~G_b=G_a$,
         mapping ~$x \in G_b$ ~to ~$ax \in G_a$,

~~~~~~ if ~$aRb$~ via isomorphism ~~$G_b.a=G_a$,
         mapping ~$x \in G_b$ ~to ~$xa \in ~~G_a$,

~~~~~~ if ~$aDb$~ via isomorphism ~~$a~G_b~a=G_a$,
         mapping ~$x \in G_b$ ~to ~$axa \in G_a$.
\end{lem}
\begin{proof}
(a) ~Let periodic element $x$ generate invariant e with period $p$, so ~$x^p=e$.
Then clearly the inverse of $x$ with respect to $e$ is $x^{p-1}$.
Define $x^{0}=e$~ for consistency in case $p$=1 ($x=e$), and denote the inverse
of $x$~ by ~$x^{-1}$. If $y$ is another periodic element generating $e$, with inverse
$y^{-1}$, then ~$xy$~ has inverse $(xy)^{-1} = y^{-1}.x^{-1}$ ~since~ $xy.(xy)^{-1}=
x.y.y^{-1}.x^{-1}= x.e.x^{-1}= x.x^{-1}= e$, and similarly $(xy)^{-1}.xy=e$.
It follows that $xy$ generates the same invariant as $x$ and $y$, so closure
holds. Inverses are unique, because if $x$ has two inverses $x_1$ and $x_2$ then
~$x_1= x_1.e= x_1.(x.x_2)= (x_1.x).x_2= e.x_2= x_2$. So all periodic elements
generating the same invariant form a group.

(b) ~Let $a,b$ be two right equivalent invariants ~$aRb$~ so ~$ab=b$~ and
~$ba=a$, then right composition of $G_a$ ~with $b$ is a morphism from $G_a$
onto $G_b$, meaning ~$G_b$ ~is an image of ~$G_a$, ~denoted ~$G_b | G_a$~
(divisor relation). This follows, because $a$ is identity for each $y$ in
$G_a: ay=ya=y$, while for each ~$x,y \in G_a: xb.yb= xb.ayb= x.ba.yb= x.a.yb=
xy.b$~ (*), where we used $ba=a$. In other words: the image of a composition
of elements is the composition of their images.

We need $ab=b$~ to show that $xb \in G_b$, in fact ~$xb$~ generates $b$ upon
iteration. This is seen by replacing $y$ in (*) with $x$, then~ $(xb)^2=
(x^2)b$, ~and in general $(xb)^i= (x^i)b$. ~Let ~$p$~ be the period of
~$x \in G_a$, ~so~ $x^p=a$, ~then ~$(xb)^p= (x^p)b= ab= b \in G_b$.

So if ~$ab=b$~ and ~$ba=a$, hence ~$a$ and $b$ are right-copiers for each other,
forming right equivalent invariants~ $aRb$, ~then right composition of ~$G_a$
with $b$ yields image ~$G_b$. Similarly, right composition of ~$G_b$~ with~
$a$~ yields image ~$G_a$. Consequently right equivalent invariants $aRb$ have
mutually ordered groups ~$G_b| G_a$~ and ~$G_a | G_b$, so they are
isomorphic:~ $G_a \cong G_b$.

Using left composition by $a$ and $b$ respectively, it follows that also left
equivalent invariants have isomorphic groups. And finally, by transitivity,
diagonal equivalent invariants have isomorphic groups as well. In that case
~$aDb$~ with (fig.3b) ~$aLba, ~~baLb$, ~and ~$a~G_b~a= a~G_{ba}= G_a$.
The diagonal case covers the other two cases of direct equivalence.
\end{proof}

{\Large\bf Conclusion}

Combining all results yields:

\begin{thm} % thm 4.1
The following conditions on a finite semigroup S are equivalent:

(a) ~$S$ is simple ~~(has no proper ideal).

(b) ~$S$ is of constant rank.

(c) ~$S$ is a disjoint union of isomorphic groups, forming
  image $L \times R$ under set-product.

(d) ~for invariants $a,b \in S: ~G_a=aSa$~ ~and~ $a~G_b~a=G_a$

(e) ~$S$ is a direct product $L \times R \times G$ ~or~ semi-direct
product $(L \times R) * G$\\
\hspace*{1cm} ~of a left- and a right-copy semigroup with a group.
\end{thm}
\begin{proof}
(a) $\Rightarrow$ (b) : Corollary 2.1 and lemma 1.1d.

(b) $\Rightarrow$ (c) : Each element $x$ of a constant rank semigroup $S$ is
periodic (cor. 4.1). Hence $S$ is a union of as many maximal subgroups as there are
invariants, being the subgroup identities (lemma 4.2a). The subgroups are disjoint
because no element can generate two invariants. Constant rank implies that no
two invariants are ordered (cor. 3.1), hence they are pairwise equivalent and form
a direct product $L \times R$ (theorem 3.1).

(c) $\Rightarrow$ (d) :
Consider an invariant ~$a$~ and elements of form ~$aSa = \{axa, ~x \in S\}$.
Let the invariant generated by ~$axa$~ be ~$c=(axa)^p$ ~with period $p$.
Since $c$ begins and ends with invariant $a$, we have ~$ac=ca=c$, meaning
~$a \geq c$, and in fact ~$a=c$, since no strict ordering occurs in a constant
rank semigroup. Hence ~$(axa)^p=a$, in other words $axa$ generates invariant $a$
for each $x$, and is thus in ~$G_a$. So for each $x$ in constant rank semigroup
$S$, ~$axa$~ is in the max-subgroup containing $a$, denoted as ~$aSa=G_a$.

If ~$a,b$~ are two equivalent invariants, with maximal subgroups $G_a$
and $G_b$, then the group isomorphism is ~$a.G_b.a=G_a$~ with ~$axa=y$,
independent of whether it is a left-, a right- or a diagonal equivalence
(lemma 4.2b), the last case covers the first two.

(d) $\Rightarrow$ (e) :
Constant rank semigroup $S$ contains as many disjoint isomorphic groups $G$ as
there are invariants. These groups form a direct product image $L \times R$
under set product (c). If the two congruences $\alpha= \{x \equiv y$
~for~ $x,y$ in the same max-subgroup\} ~and~ $\gamma=\{x \equiv y$
if $axa=y$ for some invariant $a$\} (lemma 4.2b) are orthogonal, with images
$S/ \alpha= L \times R$ and $S/ \gamma=G$, then direct product structure
$L \times R \times G$ follows. And if the product of two invariants is not
invariant then $L \times R$ is not a subsemigroup, and $G$ not an image
of $S$, yielding semi-direct product $(L \times R)*G$.

(e) $\Rightarrow$ (a) : The direct product of simple semigroups is also a
simple semigroup [2, p83, example 8]. Since ~$L, ~R$ and $G$ are simple, so
is their direct product. Although $L \times R$ is an image of $S$, it is not
necessarily a subsemigroup, in which case $G$ is not an image of $S$, with
a coupling from ($L \times R$) to $G$, corresponding to a semi-direct product.
In either case, the composition of $S$ from simple semigroups $L \times R$ and
$G$ yields transforms of equal rank, so $S$ has no proper ideal, thus is simple.
\end{proof}

Any set $A$ of state transforms that generate a constant rank closure, is
a constant rank state machine $M(A,Q)$. As shown, in general the closure
~$S= A^+/Q = (L \times R) ~*~ G$.\\
It is readily verified that $Lm$ has $m$ genetators and $m+1$ states (see $L_2$,
fig.1) with the function of an $m$-branch; $Rn$ has $n$ generators and $n$ states
with an $n$-reset function, while group $G$ has a permutation machine as generator
with $k \leq |G|$ states. Then $M$ is represented over $m+1+n+k$ states since
~$L,~R,~G$~ are 'relative prime' (have pairwise no common image, not proven here),
and we have:

\begin{cor}:\\  % cor 4.2
 A general {\bf constant rank state machine} $M$ has a semi-direct product closure
 $(L \times R ) ~* G$.\\
It is the composition of machines with closures $L, ~R, ~G$ respectively:\\
  \hspace*{1cm} a {\bf branch} machine, a {\bf reset} machine and
                a {\bf permutation} machine. \end{cor}

{\Large\bf Further research}

The decreasing-rank basic types of machines (fig.1): monotone iterative type $U$,
and combinational logic type $H$ (for instance embedding a lower semi-lattice in
a boolean lattice), still need to be included, in order to obtain a general
structure theory of State Machines. Of course, input and output logic functions
should be taken into consideration as well [3] to yield an efficient overall
logic design.

In essence, {\it associative algebra} and the {\it theory of finite semigroups}
[2] need to be translated to {\it state machine language}, and applied to
sequential logic synthesis, similar to the application of boolean algebra to
the design of combinational logic circuits. This has been tried before, but
with little practical impact, for the following reasons.

Krohn and Rhodes [4,5] derived a prime decomposition theorem using only
permutation and reset components, restricted further to cascade coupling.
This essentially extends the known Jordan-Hoelder group decomposition theorem,
by including reset machines (set/reset flipflops in the binary case). Clearly
this is not a sufficient level of detail for practical purposes: all five basic
component types [1] should be employed for a natural and efficient decomposition.

Moreover, an non-cascade or {\it loop coupling} of some components (with a
'simple group' closure) is necessary, in order to decompose such {\it 'prime'}
permutation machines, which have no casacade decomposition -- equivalent to
their sequential colure having no proper congruence. They are very complex [1]:
the smallest simple group is $A_5$, the alternating group of all 60 even
permutations of 5 states. They are not useful as practical network components.

\newpage
{\Large\bf References}
\begin{enumerate}
\item N.F.Benschop: "On State Machine Decomposition and the Five Primitives of
   Sequential Logic", {\it Internat. Workshop on Logic Synthesis}, MCNC, USA, 1987.

\item A. Clifford, G. Preston: "{\it The Algebraic Theory of Semigroups}" Vol-I,
     Mathematical Survey no.7 (AMS) 1961.

\item N.F.Benschop: "Min-Cut Algorithm for State Coding", {\it Internat. Workshop
     on Logic Synthesis}, MCNC, Research Triangle Park, NC, USA, 1989.

\item K. Krohn, J. Rhodes: "Algebraic Theory of Machines" part I,
       Tr.Am.Math.Soc, V.116, pp.450-464 (1965).

\item A. Ginzburg: "{\it Algebraic Theory of Automata}", Academic Press,
       New York 1968.
\end{enumerate}

\end{document}